\documentclass[3p]{elsarticle}

\usepackage{graphicx,amsmath,amssymb,txfonts}
\usepackage{algorithm}
\usepackage{algorithmic}
\usepackage{subfigure,color}
\usepackage[colorlinks]{hyperref}
\usepackage{hypernat}
\usepackage{multirow}
\usepackage{booktabs}
\usepackage{threeparttable}
\usepackage{epstopdf}

\graphicspath{{Figs/}}

\makeatletter

\newcommand{\Rmnum}[1]{\expandafter\@slowromancap\romannumeral #1@}
\makeatother

\newtheorem{remark}{\textbf{Remark}}

\begin{document}
\begin{frontmatter}
\title{An efficient multigrid solver for 3D  biharmonic equation with  a discretization by 25-point difference scheme}
 \author[csu]{Kejia Pan}
\author[tongji]{Dongdong He\corref{cor1}}
 \author[csu]{Runxin Ni}
 \address[csu]{School of Mathematics and Statistics, Central South University, Changsha 410083, China}
\address[tongji]{School of Aerospace Engineering and Applied Mechanics, Tongji University, Shanghai 200092, China}
  \cortext[cor1]{Corresponding authors. \\ E-mail address:dongdonghe@tongji.edu.cn(D.D. He), pankejia@hotmail.com(K.J. Pan)}

\begin{abstract}
 In this paper, we propose an efficient extrapolation cascadic multigrid (EXCMG) method combined with 25-point difference approximation to solve the three-dimensional biharmonic equation. First, through applying Richardson extrapolation and quadratic interpolation on numerical solutions on current and previous grids, a third-order approximation to the finite difference solution can be obtained and used as the iterative initial guess on the next finer grid. Then we adopt the bi-conjugate gradient (Bi-CG) method to solve the large linear system resulting from the 25-point difference approximation. In addition, an extrapolation method based on midpoint extrapolation formula is used to achieve higher-order accuracy on the entire finest grid. Finally, some numerical experiments are performed to show that the EXCMG method is an efficient solver for the 3D biharmonic equation.
\end{abstract}


\end{frontmatter}

\section{Introduction}
In this paper, we consider the following three-dimensional (3D) biharmonic equation
\begin{equation}\label{bvp}
\Delta^2 u(x,y,z)=f(x,y,z), \quad  (x,y,z) \in \Omega=(0,1)^3,
\end{equation}
with Dirichlet boundary conditions of first kind
\begin{equation}\label{BC}
u(x,y,z)=g_1(x,y,z),\quad \frac{\partial u}{\partial {n}}=g_2(x,y,z), \quad (x,y,z)\in \partial\Omega,
\end{equation}
or Dirichlet boundary conditions of second kind
\begin{equation}\label{BC2}
u(x,y,z)=g_1(x,y,z),\quad \frac{\partial^2 u}{\partial {n^2}}=g_2(x,y,z), \quad (x,y,z)\in \partial\Omega.
\end{equation}
The biharmonic operator $\Delta^2$ in three-dimensional (3D) Cartesian coordinates can be written as
\begin{equation}\label{bvp2}
\Delta^2 u(x,y,z)=\frac{\partial^4 u}{\partial x^4}+\frac{\partial^4 u}{\partial y^4}+\frac{\partial^4 u}{\partial z^4}+2\frac{\partial^4 u}{\partial x^2y^2}+2\frac{\partial^4 u}{\partial x^2z^2}+2\frac{\partial^4 u}{\partial y^2z^2}.
\end{equation}
And the two dimensional (2D) version of Eq. (\ref{bvp}) is
\begin{equation}\label{bvp2d}
  \frac{\partial^4 u}{\partial x^4}+\frac{\partial^4 u}{\partial y^4}+2\frac{\partial^4 u}{\partial x^2y^2} = f(x,y).
\end{equation}

The biharmonic equation is a fourth-order partial differential equation which arises in areas of continuum mechanics, including linear elasticity theory, phase-field models and Stokes flows. Due to the significance of the biharmonic equation, a large number of numerical methods
for solving the biharmonic equations have been proposed~\cite{Gupta19751,Gupta19752,Gupta1979,Altas1998,Bauer1972,Buzbee1974,1967differences,Altas2002,Dehghan2006,gumerov2006,FinitePointset,2017LegendreGalerkin,Conte1960,2018Schwarz,2013PGM,2014finite element}.  Most of these works  focus on two-dimensional case. There has been very little work devoted to solving the 3D biharmonic equations. The main reason is that 3D problems require large computational power and memory storage~\cite{Altas2002,Dehghan2006}.

Various methods for the numerical solutions of the biharmonic equations have been considered in the literature. A popular technique is to split $\Delta^2 u =f$ into two coupled Poisson equations for $u$ and $v$: $\Delta u=v, \Delta v=f$, each equation can be solved by using fast Poisson solvers. The coupled method has been widely used by many authors~\cite{Altas1998,Gupta19751,Gupta19752}. As it is mentioned in~\cite{Altas1998,Gupta19751,Gupta19752}, the main difficulty for the coupled (splitting or mixed) method is that the boundary conditions for the newly introduced variable $v$ are undefined and needs to be approximated accurately, and the computational results strongly depends on the choice of the approximation of missing boundary values for $v$.

 Another conventional approach for solving the 3D biharmonic equations is to directly discretize Eq. (\ref{bvp}) on a uniform grid using a 25-point computational stencil with truncation error of order $h^2$, which is derived by Ribeiro Dos Santos~\cite{1967differences} in 1967. This conventional 25-point difference approximation connects the value of $u$ at grid $(x_i, y_j, z_k)$ in terms of 24 neighboring values in a $5\times 5\times 5$ cube. Thus, this direct method need to be modified at grid points near the boundaries. As mentioned in~\cite{Gupta1979,Conte1960,Altas2002}, there are serious computational difficulties with solution of the linear systems obtained by the 13-point discretization of the 2D biharmonic equation and the 25-point discretization of 3D biharmonic equation. Dehghan and Mohebbi~\cite{Dehghan2006} also pointed that this direct method can only be used for moderate values of grid width $h$ and the well-known iterative methods such Jacobi or Gauss-Seidel either converge very slowly or diverge.

The combined compact difference method is another popular method for solving the biharmonic equation~\cite{Dehghan2006,Altas2002}. For example,  Altas et al.~\cite{Altas2002} proposed a fourth-order, combined compact formulation,  where The unknown solution and its first derivatives are carried as unkonws at grid point and computed simultaneously, for the 3D biharmonci equation with Dirichlet boundary conditions of first kind.
In 2006, Dehghan et al.~\cite{Dehghan2006} proposed two combined compact difference schemes for solve 3D biharmonic equation with Dirichlet boundary conditions of second kind, which use the known solution and its second derivatives as unknowns. In these combined compact difference methods, there is no need to modify the difference scheme at grid points near the boundaries, and the given Dirichlet boundary conditions are exactly satisfied and no approximations need to be carried out at the boundaries, in contrary to the coupled method. However, these combined compact difference methods introduce extra amount of computation, and the classical iterations for solving the resulting linear system suffer from slow convergence. Multigrid methods give good results in~\cite{Dehghan2006} and~\cite{Altas2002}. However, numerical results in~\cite{Dehghan2006} and~\cite{Altas2002} are reported only up to $32\times 32\times 32$ and $64 \times 64 \times 64$ grids, respectively. To the best of our knowledge, there is no numerical results for solving the 3D  biharmonic equations with large-scale discretized meshes.


In this paper, we propose an efficient extrapolation cascadic multigrid method based on the conventional 25-point approximation to solve 3D biharmonic equations with both first and second boundary conditions. In our method, the conventional 25-point difference scheme is used to approximate the 3D biharmonic equation (\ref{bvp}). In order to overcome the serious computational difficulties with solution of the resulting linear system, by combining Richardson extrapolation and quadratic interpolation on numerical solutions on current and previous grids, we obtain quite good initial guess of the iterative solution on the next finer grid, and then adopt the bi-conjugate gradient (Bi-CG) method to solve the large linear system efficiently. Our method has been used to solve 3D biharmonic problems with more than 135 million unknowns with only several iterations.

The rest of the paper is organized as follows: Section 2 presents the 25-point difference approximation for the 3D biharmonic equation and its modification of the difference scheme at grid points near boundaries. Section 3 reviews the classical V-cycle and W-cycle multigrid methods. In Section 4, we present a new EXCMG method to solve the linear three-dimensional biharmonic equation (\ref{bvp}). Section 5 describes the Bi-CG solver in our new EXCMG method. Section 6 provides the numerical results to demonstrate the high efficiency and accuracy of the proposed method, and conclusions are given in the final section.

\section{Second-order Finite Difference Discretization}\label{sec2}

We consider a cubic domain $\Omega=[0,1]\times[0, 1]\times[0, 1]$. Let $N = 1/h$ be the numbers of uniform
intervals along all the $x$, $y$ and $z$ directions. We discretize the domain with unequal meshsizes  $h=1/N$ in all $x, y$ and $z$ coordinate directions. The grid points are ($x_i,y_j, z_k$), with $x_i = ih, y_j = jh$ and $z_k = kh, i, j,k = 0,1,\cdots ,N$. The quantity $u_{i,j,k}$ represents the numerical solution at ($x_i,y_j, z_k$).

Then the value on the boundary points $u_{i,j,k}$ can be evaluated directly  from the Dirichlet boundary condition.
For internal grid points ($i=2,\cdots,N-2, j=2,\cdots,N-2, k=2,\cdots,N-2$), the 25-point  second-order difference scheme  for 3D biharmonic equation was derived~\cite{1967differences,Altas2002}:
\begin{align}\label{method1}
&42u_{i,j,k}-12(u_{i-1,j,k}+u_{i+1,j,k}+u_{i,j-1,k}+u_{i,j+1,k}+u_{i,j,k-1}+u_{i,j,k+1})\nonumber\\
&+u_{i-2,j,k}+u_{i+2,j,k}+u_{i,j-2,k}+u_{i,j+2,k}+u_{i,j,k-2}+u_{i,j,k+2}+\nonumber\\
&+2(u_{i-1,j-1,k}+u_{i-1,j+1,k}+u_{i+1,j-1,k}+u_{i+1,j+1,k}+u_{i-1,j,k-1}+u_{i+1,j,k-1}+u_{i,j-1,k-1}+u_{i,j+1,k-1}\nonumber\\
&+u_{i-1,j,k+1}+u_{i+1,j,k+1}+u_{i,j-1,k+1}+u_{i,j+1,k+1})=h^4f_{i,j,k}.
\end{align}

Note that $u_{i,j,k}$ is connected to grid points two grids away in each direction from the point $(x_i,y_j,z_k)$. Thus, the above difference formulation (\ref{method1}) for the grid points near the domain boundary $\partial \Omega$ involves at least one value of point outside the domain, and these points outside the domain are fictitious  points which need to be replaced by the internal points through the boundary condition. These could be done for both first and second kind of boundary conditions.

For the first kind of boundary condition. For example,  for $i=1$, the point ($x_{i-2},y_{j},z_{k}$) lies outside the computational domain, and the value on the fictitious   point ($x_{-1},y_{j},z_{k}$) can be obtained through the following central difference formula called the reflection formulas~\cite{Bauer1972, Buzbee1974}:
\begin{align}\label{reflection1}
\frac{u_{1,j,k}-u_{-1,j,k}}{2h}=\left(\frac{\partial u}{\partial x}\right)_{0,j,k},
\end{align}
where $\left(\frac{\partial u}{\partial x}\right)_{0,j,k}$ can be obtained from the boundary condition (\ref{BC}) and $u_{-1,j,k}$ is given by
\begin{align}\label{reflection2}
u_{-1,j,k}=u_{1,j,k}-2h\left(\frac{\partial u}{\partial x}\right)_{0,j,k}.
\end{align}

For the second kind of boundary condition. For example,  for $i=1$, the point ($x_{i-2},y_{j},z_{k}$) lies outside the computational domain, and the value on the fictitious   point ($x_{-1},y_{j},z_{k}$) can  also be obtained through the following central difference formula called the reflection formulas:
\begin{align}\label{reflection3}
\frac{u_{1,j,k}-2u_{0,j,k}+u_{-1,j,k}}{h^2}=\left(\frac{\partial^2 u}{\partial x^2}\right)_{0,j,k},
\end{align}
where $u_{0,j,k}$ and $\left(\frac{\partial^2 u}{\partial x^2}\right)_{0,j,k}$ can be obtained from the boundary condition (\ref{BC}) and $u_{-1,j,k}$ is given by
\begin{align}\label{reflection4}
u_{-1,j,k}=-u_{1,j,k}+2u_{0,j,k}+h^2\left(\frac{\partial^2 u}{\partial x^2}\right)_{0,j,k}.
\end{align}

We use $u_h$ and $u_{h/2}$ to represent the finite difference solutions of equation (\ref{bvp}) with mesh sizes $h$ and $\frac{h}{2}$ respectively. Afterward, a matrix form, which express the finite difference scheme (\ref{method1}) and an equation set including formulas of the grid points near the boundary, can be obtained as below:
\begin{equation}\label{sec}
A_h u_h=f_h,
\end{equation}
Where $A_h$ is not a symmetry positive definite matrix, and the right hand-side vector of (\ref{method1}) and an equation set including the formulas of the grid points near the boundary are expressed by $f_h$.

Note that the discretization equations for grid points that away from the boundary and that near the boundary are different, one must distinguish all possible cases. Although there are a little bit troublesome to treat all cases (there are totally 27 cases with 27 different equations), by moving the known boundary values into the right hand-side of the system, it is convenient to solve these equations which only involves unknown on the grid points.


\section{Classical Multigrid Method}\label{sec3}

Since the 1970s', many scholars have done researches on the classical multigrid method. Through deep researches on it for about fifty years, the classical multigrid method gradually forms its own comprehensive system. Including the interpolation, restriction and iteration, the classical multigrid method starts from the fine grid, goes to coarse grid and then returns to the fine grid. The classical multigrid methods contain V-cycle and W-cycle.

\begin{figure}[!tbp]
  \centering
  \includegraphics[width=5in]{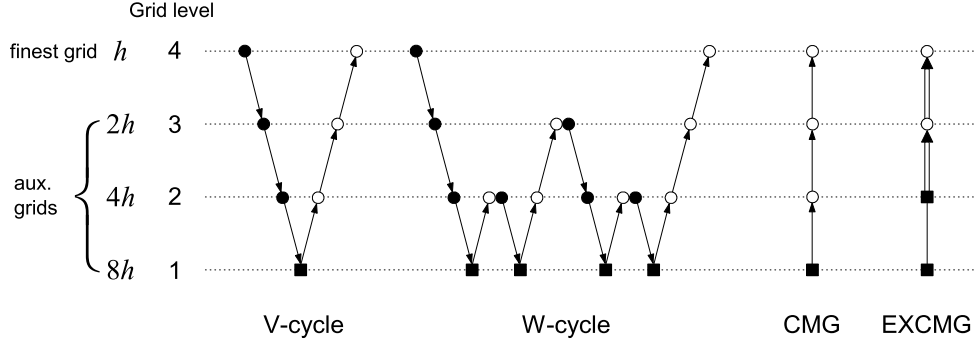}\\
  \caption{The four level structure of the V- and W-cycles, CMG and EXCMG methods. In the diagram, $\bullet$ denotes pre-smoothing steps, $\circ$ denotes post-smoothing steps, $\uparrow$ denotes prolongation, $\downarrow$ denotes restriction, $\Uparrow$ denotes extrapolation and quadratic interpolation, and $\blacksquare$ denotes direct solver.}
  \label{VW}
  \end{figure}

The classical multigrid method is introduced in detail with several steps. First, the specific smoother is used to smooth the current approximation on the fine grid. To obtain more oscillatory error components, we compute the residual and transfer it to the coarser grid with restriction. Next, we solve the residual equation on the coarser grid with the application of the number ($\gamma$) of cycles. From the fine grid to the coarsest grid and back to the fine grid is called a cycle. Then, we acquire the improved approximation on the fine grid by interpolating the correction back to the fine grid. Finally, we smooth the obtained approximation on the fine grid with the smoother again. If $\gamma$=1, call it V-cycle. And if $\gamma$=2, call it W-cycle. We take the four-level structures of V-cycle and W-cycle in Fig.\ref{VW} for instances to illustrate that.

\begin{remark}
When the $\gamma$-cycle is performed on the coarsest grid, direct solver is used to solve the residual equation.
\end{remark}

 \section{Extrapolation Cascadic Multigrid Methods}\label{sec4}
It is an important issue to find approaches to solving the linear equation with enormous unknowns, which is obtained by FE and FD discretizations. Therefore, many authors paid great attention on it and presented multigrid methods including the MG method, the CMG method and the EXCMG method. The MG method has had a nearly integrated system through many scholars' hard work in the past several decades. However, its algorithm is complex. Then the CMG method proposed by Deuflhard and Bornemann in~\cite{1996CMG} only use the interpolation and iteration so that its algorithm which is easy to operate is appealing. Furthermore, in 2008, the EXCMG method was proposed by Chen et al.~\cite{2008EXCMG} and the cores of it are Richardson extrapolation and quadric interpolation. Compared with the CMG method, the EXCMG method provides a much better initial guess for the iteration solution on the next finer grid. In this section, we propose a new EXCMG method combined with the second-order incompact FD discretization for solving the linear three-dimensional biharmonic equation.

\subsection{Description of the EXCMG Algorithm}

In Algorithm \ref{alg:EXCMG}, \emph{H} is the size of the coarsest grid.  \emph{L}, the positive integer, denotes the total number of grids except first two embedded grids and indicates that the finest grids' size is  $\frac{H}{2^{L+1}}$. For the sizes of first two coarse grids are small, DSOLVE, a direct solver, is applied on the first two coarse grids (see line 1-2 in the Algorithm \ref{alg:EXCMG}). In addition, procedure $\textrm{EXP}_{finite}(u_{2h}, u_{4h})$ represents the third-order approximation of the FD solution $u_h$ which is obtained by Richardson extrapolation and quadratic interpolation from numerical solutions $u_{2h}$ and $u_{4h}$. Meanwhile, a selective step is presented in the Algorithm \ref{alg:EXCMG} above where $\textrm{EXP}_{true}(u_{h}, u_{2h})$ refers to a higher-order solution which is extrapolated on the finest grid with the mesh size, $h$, from two second-order numerical solutions $u_h$ and $u_{2h}$.

\begin{algorithm}[!tbp]

\caption{New EXCMG : $(u_h, \tilde{u}_h)$ $\Leftarrow$ EXCMG($A_h, f_h, L ,\epsilon$)}

\label{alg:EXCMG}

\begin{algorithmic}[1]

\STATE $u_H$ $\Leftarrow$ DSOLVE($A_H u_H=f_H$)
\STATE $u_{H/2}$ $\Leftarrow$ DSOLVE($A_{H/2} u_{H/2}=f_{H/2}$)
\STATE $h=H/2$

  \FOR {$i=1$ to $L$}
   \STATE $h = h/2$, $\quad \epsilon_i = \epsilon \cdot 10^{i-L}$
   \STATE ${w}_{h} = \textrm{EXP}_{finite}(u_{2h}, u_{4h})$,  $\quad \quad \quad \rhd$ $u_h=w_h$ is used as the initial guess  for Bi-CG solver
    \WHILE {$||A_h u_h -f_h||_2>\epsilon_i \cdot ||f_h||_2 $}
        \STATE $u_h \Leftarrow$ Bi-CG$(A_h, u_h, f_h)$
    \ENDWHILE
  \ENDFOR
  \STATE $\tilde{u}_{h} = \textrm{EXP}_{true}(u_{h}, u_{2h})$  $\quad \quad \quad \quad \quad \rhd$ $\tilde{u}_{h}$ is a higher-order approximation solution
\end{algorithmic}
\end{algorithm}

The details of the procedure of extrapolation and quadratic interpolation are introduced next subsection \ref{extra}. The difference between our new EXCMG method and existing EXCMG method are discussed below:

\begin{enumerate}[(1)]
\item Instead of applying the second-order linear FE method, a second-order incompact difference scheme is used to discretize the 3D biharmonic equation in our new EXCMG  method.

\item Rather than perform the fixed number of iterations used in the existing EXCMG method, we introduce a relative residual tolerance $\epsilon_i$ into the Bi-CG solver (see line 7 in the Algorithm \ref{alg:EXCMG}), which enables us to avoid the difficulty of determining the number of iterations at every grid level and obtain numerical solutions with desired accuracy conveniently.

\item In our new EXCMG method, we take the Bi-CG solver as smoother instead of the CG solver (see line 8 in Algorithm \ref{alg:EXCMG}). The Bi-CG is more suitable for positive definite matrix which is not symmetric compared with the CG solver.

\item Through $\textrm{EXP}_{true}(u_{h}, u_{2h})$, a higher-order extrapolated solution $\tilde{u}_{h}$ is obtained easily, which improves the accuracy of the numerical solution $u_h$ (see line 11 in Algorithm \ref{alg:EXCMG}).

\end{enumerate}

\subsection{Extrapolation and Quadratic Interpolation}\label{extra}

The Richardson extrapolation is a well-known method for producing more accurate solutions of many problems in numerical analysis. Marchuk and Shaidurov~\cite{Marchuk1983} researched the application of the Richardson extrapolation on the FD method systematically in 1983. Since then, this technique has been well demonstrated in the frame of the FE and FD methods~\cite{2010sixth, Marchuk1983, 1987Acceleration, 1989regular, 1993Spline, 2004FD, 2006FD, 2008R, 2000wave, 20103D, 20132D}.

In next three subsections, we will give the explanation for how to obtain higher-oder accuracy solution on the fine grid. Moreover, how to acquire a third-order approximation of the second-order FD method on the next finer grid is illustrated as well. Meanwhile, we can regard it as another critical application of the extrapolation method which produces good initial guesses for iterative solutions.

\subsubsection{Extrapolation for the True Solution}

For simplicity, we first consider the three-level of embedded grids \emph{$Z_i$}(\emph{i}=0, 1, 2) with mesh sizes \emph{$h_i$} = \emph{$h_0$}/$2^i$ in one dimension. In addition, let $e^i$ = $u^i$-\emph{u} be the error of the second-order incompact FD solution $u^i$ with mesh size $h_i$. We make an assumption that the error at the node has the following form:

\begin{equation}\label{secf}
  e^i(x_k)=A(x_k)h_i^2+O(h_i^4),
\end{equation}

where \emph{A}(\emph{x}) is a properly smooth function independent of $h_i$. We will verify the error expansion (\ref{secf}) by numerical results in Sect. 5.

Through the equation (\ref{secf}), the Richardson extrapolation formula at the coarse grid point is obtained

\begin{equation}\label{t1}
  \tilde{u}_k^1 := \frac{4 u^1_k - u^0_k}{3} = u(x_k) + O(h_0^4),\ \ k=j,j+1.
\end{equation}

Then, a midpoint extrapolation formula is obtained by linear interpolation

\begin{equation}\label{chenlin}
  \tilde{u}_{j+1/2}^1 := u_{j+1/2}^1 + \frac{1}{6}(u^1_j - u^0_j + u^1_{j+1} - u^0_{j+1}) = u(x_{j+1/2}) + O(h_0^4),
\end{equation}

whose accuracy is fourth-order at fine grid points.

From equation (\ref{secf}), it is easy to obtain

\begin{equation}\label{aaa}
  A(x_k) = \frac{4}{3h_0^2}(u_k^0-u_k^1) + O(h_0^2),\ \ k=j,j+1.
\end{equation}

Through the error estimate of the linear interpolation

\begin{equation}\label{A_12}
  A(x_{j+1/2})=\frac{1}{2}(A(x_j)+A(x_{j+1}))+O(h_0^2).
\end{equation}

and substituting equation (\ref{aaa}) to equation (\ref{A_12}), yield

\begin{equation}\label{Axmid}
  A(x_{j+1/2}) = \frac{2}{3h_0^2}(u_j^0-u_j^1) + \frac{2}{3h_0^2}(u_{j+1}^0-u_{j+1}^1) + O(h_0^2).
\end{equation}

Since

\begin{equation}
  u_{j+1/2}^1 = u(x_{j+1/2}) + \frac{1}{4} A(x_{j+1/2})h_0^2 + O(h_0^4),
\end{equation}

The midpoint extrapolation formula (\ref{chenlin}) is obtained by the use of equation (\ref{Axmid}).

\subsubsection{Extrapolation for the FD Solution}

In this subsection, given solutions $u^0$ and $u^1$ of the second-order FD method, we will explain how to construct a third-order approximation $w^2$ of the FD solution $u^2$ by using extrapolation and interpolation methods in detail.

\begin{figure}[!tbp]
   \centering
   \scalebox{0.5}{\includegraphics{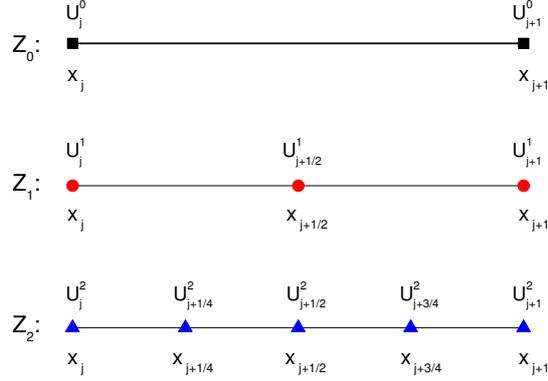}}
   \caption{Three embedded mesh in 1D.}\label{Fig:1}
\end{figure}

We divide the coarse element ($x_j$,$x_{j+1}$) into four uniform elements by adding one midpoint and two four equal points which are on the left side and right side of the midpoint. As a result, a set which contains five points and belongs to fine mesh $Z_2$ is obtained $$\Big\{x_j,x_{j+1/4},x_{j+1/2},x_{j+3/4},x_{j+1}\Big\}$$

To acquire the more accurate approximation of FD solution $u^2$, the given solutions $u^0$ and $u^1$ are combined linearly. Therefore, here assume the existence of a constant \emph{c} such that

\begin{equation}\label{secff}
  c u^1+(1-c) u^0=u^2+O(h_0^4),
\end{equation}

For obtaining the value of the constant \emph{c}, we substitute the error expansion (\ref{secf}) into  (\ref{secff}) and obtain \emph{c}=5/4. Afterward, we obtain formulas of node extrapolation at points $x_j$ and $x_{j+1}$.

\begin{equation}\label{secfff}
 w^2_k := \frac{5 u^1_k - u^0_k}{4} = u^2_k+O(h_0^4),\ \ k=j,j+1,
\end{equation}

Next, derivate the midpoint $x_{j+1/2}$'s extrapolation formula. First, use the error expansion (\ref{secf}) again and obtain the formula below

\begin{equation}\label{fsec}
u_{j+1/2}^2 = u_{j+1/2}^1 - \frac{3}{16} A(x_{j+1/2}) h_0^2 + O(h_0^4),
\end{equation}

Then substitute the (\ref{Axmid}) to (\ref{fsec}) for eliminating the unknown \emph{A}($x_{j+1/2}$) and the following extrapolation formula of $x_{j+1/2}$ is yielded

\begin{equation}\label{bfsec}
w^2_{j+1/2}:=u^1_{j+1/2}+ \frac{1}{8}(u^1_j - u^0_j + u^1_{j+1} - u^0_{j+1})={u}^2_{j+1/2} + O(h_0^4).
\end{equation}

Finally, since the values of three points $w_j^2$, $w_{j+1/2}^0$ and $w_{j+1}^0$ have been obtained, we can derivate the extrapolation formulas at the points $x_{j+1/4}^2$ and $x_{j+3/4}^2$ blow through the use of the quadratic interpolation method

\begin{align}
      w^2_{j+1/4}&:=\displaystyle\frac1{16}\big[(9u^1_{j}+12u^1_{j+1/2}-u^1_{j+1})-(3u^0_j+u^0_{j+1})\big],\label{sifen1}\\
      w^2_{j+3/4}&:=\displaystyle\frac1{16}\big[(9u^1_{j+1}+12u^1_{j+1/2}-u^1_j)-(3u^0_{j+1}+u^0_j)\big].\label{sifen2}
\end{align}

From the polynomial interpolation's theory, it is easy to demonstrate that the third-order approximation of the FD solution can be presented by formulas (\ref{efsec}) and (\ref{ffsec}). i.e.,

\begin{equation}\label{efsec}
w^2_{j+1/4}=u^2_{j+1/4} + O(h_0^3),
\end{equation}
\begin{equation}\label{ffsec}
w^2_{j+3/4}=u^2_{j+3/4} + O(h_0^3),
\end{equation}

\subsubsection{Application of Extrapolation and Quartic Interpolation on Three-Dimension}

\begin{figure}[!tbp]
  \centering
  \subfigure[$Z_0$]{
  \begin{minipage}[t]{.3\linewidth}
\includegraphics[width=0.95\textwidth]{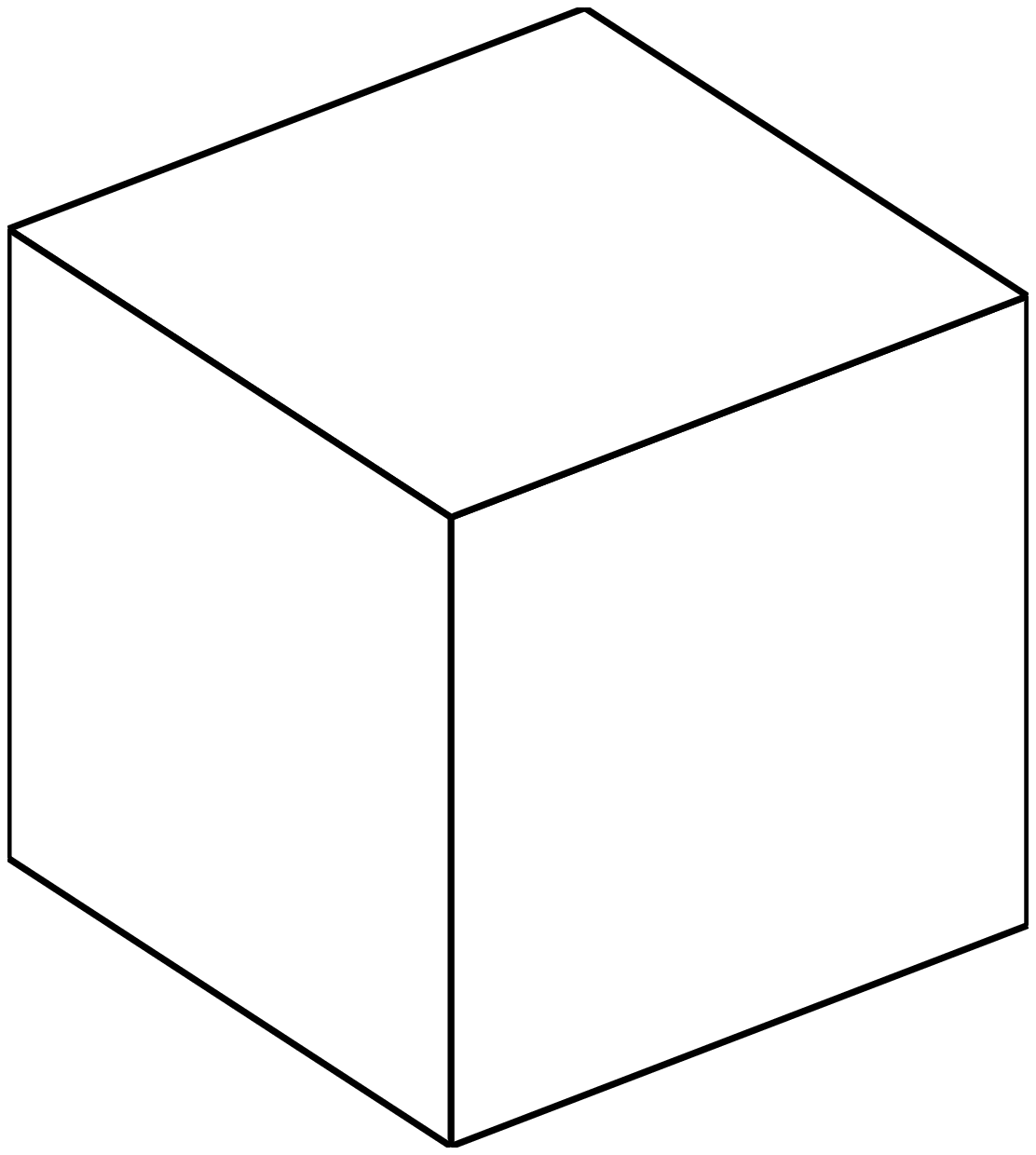}
  \end{minipage}
  }
  \subfigure[$Z_1$]{
  \begin{minipage}[t]{.3\linewidth}
\includegraphics[width=0.95\textwidth]{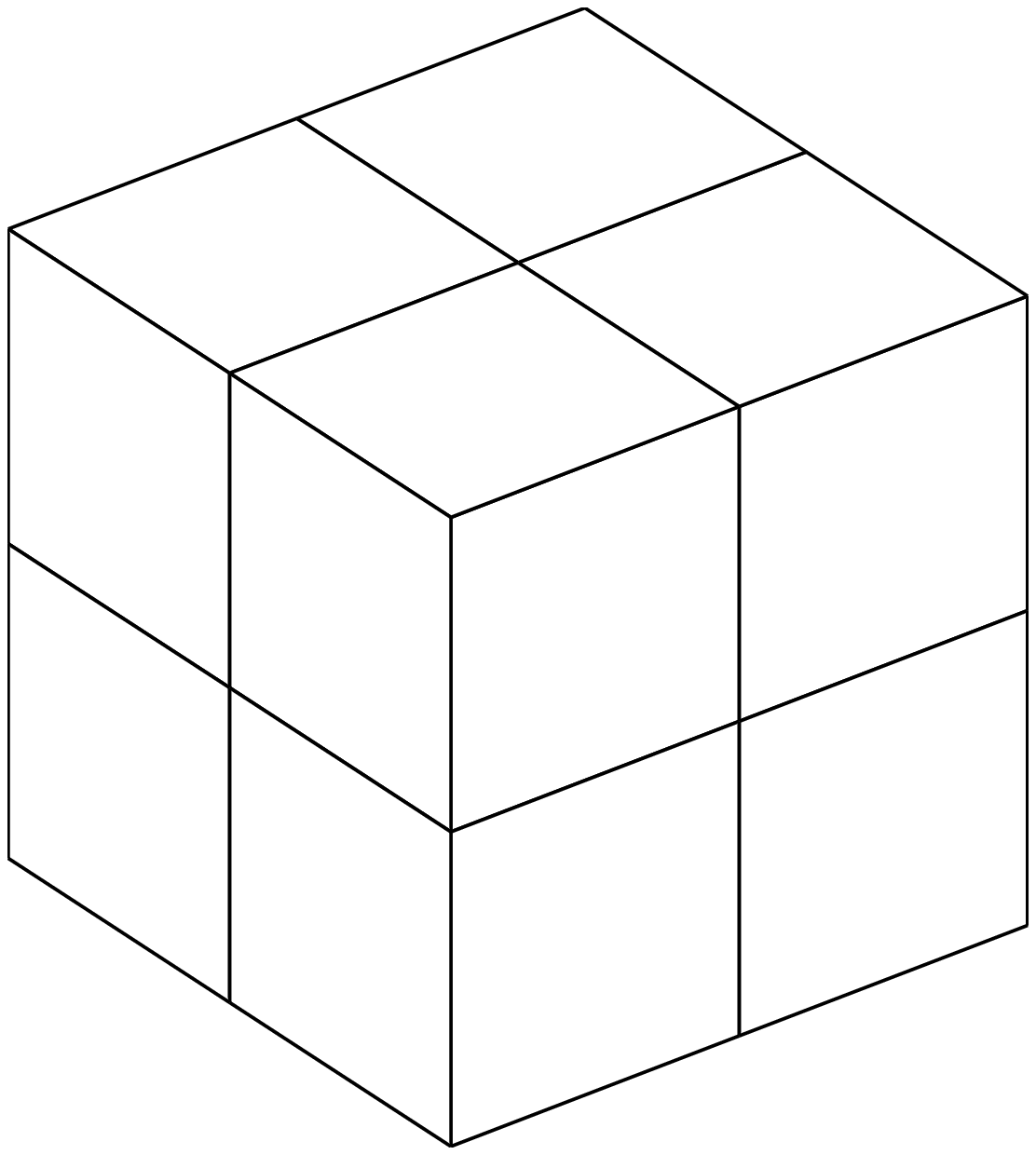}
  \end{minipage}
  }
  \subfigure[$Z_2$]{
  \begin{minipage}[t]{.3\linewidth}
\includegraphics[width=1.0\textwidth]{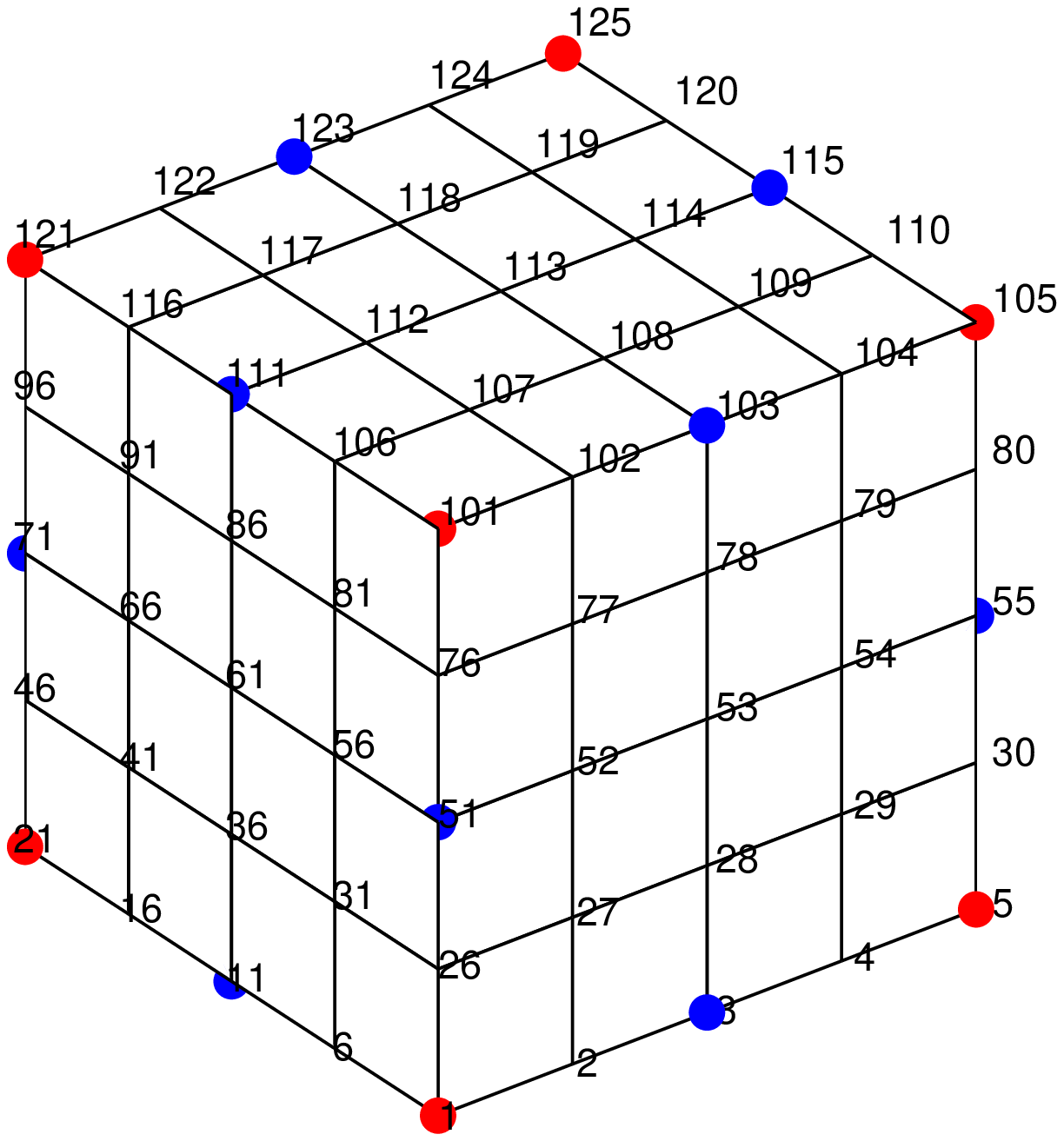}
  \end{minipage}
  }
  \caption{Three embedded hexahedral mesh}\label{Fig3D}
\end{figure}

In this subsection, how to acquire an accurate third-order approximation $w^2$ of the FD solution $u^2$ will be explained for embedded cubic mesh as shown in Fig.\ref{Fig3D}. The specific steps of the construction of the third-order approximation $u^2$ are illustrated below:\\
\textbf{Step 1}: Corner nodes (1, 5, 21, 25, 101, 105, 121, 125): Use the extrapolation formula (\ref{secfff}) to obtain the approximations of the values of 8 corner nodes on interpolation cells.\\
\textbf{Step 2}: Midpoints of edges (3, 11, 15, 23, 51, 55, 71, 75, 103, 111, 115, 123): Use the midpoint extrapolation formula (\ref{bfsec}) to obtain the approximations of the values of 19 midpoints of edges on interpolation cells.\\
\textbf{Step 3}: Centers of faces (13, 53, 65, 61, 77, 113): View the center of each face as the midpoint of the two face diagonals on interpolation cells. To obtain the approximations of the values of them, use the midpoint extrapolation formula (\ref{bfsec}) obtaining two approximations, calculate the arithmetic mean of the two obtained approximations and treat it as the approximate value at the center of each face.\\
\textbf{Step 4}: Center of the hexahedral element (63): View the center of the hexahedral element as the midpoint of the four space diagonals on interpolation cells. To obtain the approximation of the value of it, use the midpoint extrapolation formula (\ref{bfsec}) obtaining four approximations, calculate the arithmetic mean of the four obtained approximations and treat it as the approximate value at the hexahedral element.\\
\textbf{Step 5}: Other 98 fine grid points:  The approximations of remaining 98 ($5^3$ - 27) grid points can be obtained by using tri-quadratic interpolation with the known 27-node (8 corner nodes, 12 midpoints of edges, 6 centers of faces and 1 center of the hexahedral element) values.

The tri-quartic interpolation formula at natural coordinates ($\xi, \eta, \zeta$) is defined as
  \begin{equation}\label{inter}
      w^2(\xi, \eta, \zeta)=\sum_i N_i(\xi, \eta, \zeta)w_i^2,
    \end{equation}
    where the shape functions $N_i$ can be written below:
\begin{equation}\label{ininter}
    {N_i}(\xi ,\eta ,\zeta) = {l_i^2}(\xi){l_i^2}(\eta){l_i^2}(\zeta)
  \end{equation}
where $l_i^2$(\emph{x}) (0 $\leq$ \emph{i} $\leq$ 2) are the Lagrange basis polynomials of degree 2 which are defined as
\begin{equation}\label{iinter}
    {l_i^2}(\xi)=\prod_{k=0,k\not=i}^2\frac{\xi-\xi_k}{\xi_i-\xi_k}
  \end{equation}
and ($\xi$, $\eta$, $\zeta$) is the natural coordinate of node \emph{i} (1 $\leq$ \emph{i} $\leq$ 27).

\section{Bi-Conjugate Gradient Method}\label{sec5}

The Bi-Conjugate Gradient (Bi-CG) method is an algorithm which is focus on solving linear equation systems
\begin{equation}\label{1inter}
   \emph{Ax} = \emph{b}
  \end{equation}
Compared with the Conjugate Gradient (CG) method which needs matrix \emph{A} to be self-joint, the Bi-CG method does not require matrix \emph{A} self-joint but require it to multiply conjugate transpose $A^{*}$. In addition, the Bi-CG method replaces the residual's orthogonal sequence in the CG method with two sequences which are mutually orthogonal. In the Bi-CG method, the residual $r_j$ is orthogonal with a set of vectors $\hat{r}_0$, $\hat{r}_1$ ... $\hat{r}_{j-1}$ and $\hat{r}_j$ is also orthogonal with $r_0$, $r_1$ ... $r_{j-1}$. These relationships can be achieved by two three-term recurrence relations of vectors $\{$$\hat{r}_j$$\}$ and $\{$$r_j$$\}$. Meanwhile, the Bi-CG method terminates within at most \emph{n} steps when \emph{A} is an \emph{n} by \emph{n} matrix. The preconditioned version and the unpreconditioned version of algorithms of the Bi-CG method are described as follows:\\

In the Bi-CG Method with the preconditioner algorithm below, $x_k^{*}$ is adjoint, $\bar{\alpha}$ is the complex conjugate and the calculated $r_k$ and  $r_k^{*}$  satisfy the following equations respectively
\begin{align}
   r_k = \emph{b}-\emph{A$x_k$}\label{1}\\
   r_k^{*} = \emph{$b^{*}$}-\emph{$x_k^{*}$A}\label{2}
  \end{align}

\begin{algorithm}
\caption{Bi-CG Method with the Preconditioner}

\begin{algorithmic}[2]

\STATE $x_0$ is an initial guess
\STATE Choose two other vectors $x_0^{*}$ and $b^{*}$ and a preconditioner \emph{M}
\STATE $r_0$=\emph{b}-\emph{A$x_0$}
\STATE $r_0^{*}$=\emph{$b^{*}$}-\emph{$x_0^{*}$A}, such that ($r_0$, $r_0^{*})$ $\not=$ 0
\STATE $p_0$=\emph{$M^{-1}$}$r_0$
\STATE $p_0^{*}$=$r_0^{*}$\emph{$M^{-1}$}

  \FOR {$k=0,1,...$}
\STATE $\alpha_k$=$\frac{r_k^{*}\emph{$M^{-1}$}r_k}{p_k^{*}\emph{A}p_k}$
\STATE $x_{k+1}$=$x_k$+$\alpha_k$$p_k$
\STATE $x_{k+1}^{*}$=$x_k^{*}$+$\bar{\alpha_k}$$p_k^{*}$
\STATE $r_{k+1}$=$r_k$-$\alpha_k$\emph{A}$p_k$
\STATE $r_{k+1}^{*}$=$r_k^{*}$-$\bar{\alpha_k}$$p_k^{*}$\emph{A}
\STATE $\beta_k$=$\frac{r_{k+1}^{*}\emph{$M^{-1}$}r_{k+1}}{r_k^{*}\emph{$M^{-1}$}r_k}$
\STATE $p_{k+1}$=\emph{$M^{-1}$}$r_{k+1}$+$\beta_k$$p_k$
\STATE $p_{k+1}^{*}$=$r_{k+1}^{*}$\emph{$M^{-1}$}+$\bar{\beta_k}$$p_k^{*}$
    \ENDFOR
\end{algorithmic}
\end{algorithm}

%
%
%

In this paper, we adopt the Bi-CG method with  the preconditioner as the iteration solver in our new EXCMG method.

\section{Numerical Experiments}\label{sec6}

{\bf Test Problem 1.}
The exact solution of the test problem 1 introduced in~\cite{Numer} can be written as
\begin{equation}\label{test1}
u(x,y,z)= (1-\cos(2{\pi}x))(1-\cos(2{\pi}y))(1-\cos(2{\pi}z)).
\end{equation}
Applying the biharmonic operator on the exact solution, we can obtain the forcing term \emph{f}(\emph{x}, \emph{y}, \emph{z}) as follows:
\begin{equation}\label{force1}
\begin{split}
\emph{f}(\emph{x}, \emph{y}, \emph{z})=-16{\pi}^4(\cos(2{\pi}x)-4\cos(2{\pi}x)\cos(2{\pi}z)-4\cos(2{\pi}x)
\cos(2{\pi}y)\\
+9\cos(2{\pi}x)\cos(2{\pi}y)\cos(2{\pi}z)+\cos(2{\pi}y)-4\cos(2{\pi}y)\cos(2{\pi}z)+\cos(2{\pi}z))
\end{split}
\end{equation}

Obtain first boundary data from the exact solution while obtaining second boundary data \emph{$f_2$} by taking partial derivative for the exact solution.

Results listed on the table \ref{table1} of numerical experiments are performed with  EXCMG$_{bi-cg}$, using 3.6 thousand unknowns on the coarsest grid 32$\times$32$\times$32 and more than 135 million unknowns on the finest grid 512$\times$512$\times$512. In the table \ref{table1}, ``Iter'' denotes the number of iterations needed for the Bi-CG solver to achieve the relative residual less than the given tolerance. Additionally, the last row in the table provides the the $L^{\infty}$-error and $L^2$-error of the extrapolated solution $u_h$ on the finest grid, and the amount of computational cost of the ${EXCMG}_{bi-cg}$ method in terms of a work unit (WU) on the finest grid, which is defined as the total computation required to perform one relaxation sweep on the finest grid. We use the same notations in all tables.

From the results in the table \ref{table1}, it is clear that the numerical solution $u_h$ reaches almost full second-order accuracy while the initial guess $w_h$ is third-order approximation to numerical solution $u_h$. In addition, the extrapolated solution $\tilde{u}_{h}$ increase the numerical solution's accuracy greatly. What's more, the number of iterations is reduced significantly while the grids are finer and finer and this feature is especially important while solving large linear systems. We will introduce this feature specifically in the following text.

First, we we define the error ratio $r_h$ as
\begin{equation}\label{ratiorh}
r_h=\frac{||w_h-u_h||_{2}}{||u_h-u||_{2}}
\end{equation}
For the order of $||w_h-u_h||_{2}$ is one higher than the order of $||u_h-u||_{2}$, the error ratio $r_h$ is almost $\frac{1}{2^n}$ where n denotes the level of the grid. As the grid becomes finer, $w_h$ is much closer to $u_h$ especially on the finest grid. Therefore, when the grid is fine enough, the error of $||u_h-w_h||_{2}$ is smaller so that the number of iterations is reduced.

For test problem 1, on the finest grid 512$\times$512$\times$512, the error ratio $r_h$ is 0.028. It is obvious that the $r_h$ is so small on the finest grid that we only need to perform one iteration on the finest grid. The number of iteration is reduced significantly.


\begin{table}[!tbp]
\tabcolsep=6pt
\caption{Errors and convergence rates using EXCMG$_{bi-cg}$ for Problem 1.} \centering
\begin{threeparttable}
  \begin{tabular}{|c|c|cc|cc|cc|c|}
    \hline
    \multirow{2}{*}{Mesh}& \multirow{2}{*}{Iters} & \multicolumn{2}{c|}{$||u_h-u||_{2}$}  & \multicolumn{2}{c|}{$||u_h-u||_{\infty}$} & \multicolumn{2}{c|}{$||w_h-u_h||_{2}$}  \\
    \cline{3-8}&   &     Error &   Order  &     Error &   Order  &   Error &   Order  \\
\hline
$  32\times  32\times  32$ &   474 &$1.13(-2)$ &         &$5.16(-2)$ &       &  $5.14(-3)$ &        \\
$  64\times  64\times  64$ &   512 & $2.89(-3)$ & 1.97  &  $1.29(-2)$ & 2.00  &  $6.15(-4)$ & 3.06   \\
$ 128\times 128\times 128$ &   64 & $7.27(-4)$ & 1.99  &  $3.21(-3)$ & 2.00  &  $7.63(-5)$ & 3.01   \\
$ 256\times 256\times 256$ &   8 & $1.80(-4)$ & 2.01 & $8.01(-4)$ & 2.00  &   $9.56(-6)$ & 3.00   \\
$ 512\times 512\times 512$ &   1 & $4.25(-5)$ & 2.08   &  $1.99(-4)$ & 2.01 &   $1.19(-6)$ & 3.00    \\
\hline
  &  4.12 WU          & &  &  $1.46(-5)$ & &  $5.81(-6)$ &                           \\
\hline
\end{tabular}
\begin{tablenotes}
        \footnotesize
        \item[1] WU ({\it work unit}) is the computational cost of performing one relaxation sweep on the finest grid. Here, the $\textrm{EXCMG}_{bi-cg}$
        computation cost $= 1 + 8\times 2^{-3} + 64\times 2^{-6} + 512\times 2^{-9} + 474\times 2^{-12} \approx 4.12$.
      \end{tablenotes}
\end{threeparttable}
\label{table1}
\end{table}

~\\
{\bf Test Problem 2.}
The exact solution of the test problem 2 can be written as
\begin{equation}\label{test2}
u(x,y,z)= e^{xyz}.
\end{equation}
Applying the biharmonic operator on the exact solution, we can obtain the forcing term \emph{f}(\emph{x}, \emph{y}, \emph{z}) as follows:
\begin{equation}\label{force2}
\emph{f}(\emph{x}, \emph{y}, \emph{z})=e^{xyz}(x^{4}y^{4}+2x^{4}y^{2}z^{2}+x^{4}z^{4}+8x^{3}yz+2x^{2}y^{4}z^{2}+2x^{2}y^{2}z^{4}+4x^{2}+8xy^{3}z+8xyz^{3}+y^{4}z^{4}+4y^{2}+4z^{2})
\end{equation}

Obtain first boundary data from the exact solution while obtaining second boundary data \emph{$f_2$} by taking partial derivative for the exact solution.

Again, results listed on the table \ref{table2} of numerical experiments are performed on five level grids with 3.6 thousand unknowns on the coarsest grid 32$\times$32$\times$32 and more than 135 million unknowns on the finest grid 512$\times$512$\times$512. Moreover, from table \ref{table2}, we can see that numerical solution $u_h$ reaches almost full second-order accuracy, the initial guess $w_h$ is third-order approximation to numerical solution $u_h$, while the extrapolated solution $\tilde{u}_{h}$ increase the numerical solution's accuracy significantly. If we use $w_h$ as numerical solution on the finest grid 512$\times$512$\times$512, the error ratio $r_h$ is already 0.27. Thus only six iterations are needed to perform to achieve the expected accuracy.

\begin{table}[!tbp]
\tabcolsep=6pt
\caption{Errors and convergence rates using EXCMG$_{bi-cg}$ for Problem 2.} \centering
\begin{threeparttable}
  \begin{tabular}{|c|c|cc|cc|cc|c|}
    \hline
    \multirow{2}{*}{Mesh}& \multirow{2}{*}{Iters} & \multicolumn{2}{c|}{$||u_h-u||_{2}$}  & \multicolumn{2}{c|}{$||u_h-u||_{\infty}$} & \multicolumn{2}{c|}{$||w_h-u_h||_{2}$}  \\
    \cline{3-8}&   &     Error &   Order  &     Error &   Order  &   Error &   Order  \\
\hline
$  32\times  32\times  32$ &   259 &$8.96(-7)$ &         &$8.06(-6)$ &       &  $4.59(-6)$ &        \\
$  64\times  64\times  64$ &   470 & $2.30(-7)$ & 1.96  &  $2.06(-6)$ & 1.97  &  $5.48(-7)$ & 3.07   \\
$ 128\times 128\times 128$ &   384 & $5.80(-8)$ & 1.99  &  $5.15(-7)$ & 2.00  &  $6.66(-8)$ & 3.04   \\
$ 256\times 256\times 256$ &   48 & $1.46(-8)$ & 1.99 & $1.28(-7)$ & 2.01  &   $8.19(-9)$ & 3.02   \\
$ 512\times 512\times 512$ &   6 & $3.67(-9)$ & 1.99   &  $3.22(-8)$ & 1.99 &   $1.00(-9)$ & 3.03    \\
\hline
  &  18.98 WU          & &  &  $9.90(-9)$ & &  $2.93(-10)$ &                           \\
\hline
\end{tabular}
\begin{tablenotes}
        \footnotesize
        \item[1] WU ({\it work unit}) is the computational cost of performing one relaxation sweep on the finest grid. Here, the $\textrm{EXCMG}_{bi-cg}$
        computation cost $= 6 + 48\times 2^{-3} + 384\times 2^{-6} + 470\times 2^{-9} + 259\times 2^{-12} \approx 18.98$.
      \end{tablenotes}
\end{threeparttable}
\label{table2}
\end{table}

~\\
{\bf Test Problem 3.}
The exact solution of the test problem 3 can be written as
\begin{equation}\label{test3}
u(x,y,z)= \sinh(x)\sinh(y)\sinh(z).
\end{equation}

Applying the biharmonic operator on the exact solution, we can obtain the forcing term \emph{f}(\emph{x}, \emph{y}, \emph{z}) as follows:
\begin{equation}\label{force3}
\emph{f}(\emph{x}, \emph{y}, \emph{z})=\sinh(x)\sinh(y)\sinh(z).
\end{equation}

Obtain first boundary data from the exact solution while obtaining second boundary data \emph{$f_2$} by taking partial derivative for the exact solution.

Again, five level grids are used with 3.6 thousand unknowns on the coarsest grid 32$\times$32$\times$32 and more than 135 million unknowns on the finest grid 512$\times$512$\times$512. In addition, from table \ref{table3}, we can see that numerical solution $u_h$ reaches almost full second-order accuracy, the initial guess $w_h$ is third-order approximation to numerical solution $u_h$, while the extrapolated solution $\tilde{u}_{h}$ increase the numerical solution's accuracy greatly. On the finest grid 512$\times$512$\times$512, if using $w_h$ as numerical solution, the error ratio $r_h$ is already equal to 0.13. Thus we only need to perform six iterations to achieve the expected accuracy.

\begin{table}[!tbp]
\tabcolsep=6pt
\caption{Errors and convergence rates using EXCMG$_{bi-cg}$ for Problem 3.} \centering
\begin{threeparttable}
  \begin{tabular}{|c|c|cc|cc|cc|c|}
    \hline
    \multirow{2}{*}{Mesh}& \multirow{2}{*}{Iters} & \multicolumn{2}{c|}{$||u_h-u||_{2}$}  & \multicolumn{2}{c|}{$||u_h-u||_{\infty}$} & \multicolumn{2}{c|}{$||w_h-u_h||_{2}$}  \\
    \cline{3-8}&   &     Error &   Order  &     Error &   Order  &   Error &   Order  \\
\hline
$  32\times  32\times  32$ &   285 &$4.10(-6)$ &         &$1.75(-5)$ &       &  $9.49(-6)$ &        \\
$  64\times  64\times  64$ &   533 & $1.05(-6)$ & 1.96  &  $4.36(-6)$ & 2.00  &  $1.15(-6)$ & 3.04   \\
$ 128\times 128\times 128$ &   384 & $2.66(-7)$ & 1.98  &  $1.09(-6)$ & 2.00  &  $1.42(-7)$ & 3.02   \\
$ 256\times 256\times 256$ &   48 & $6.71(-8)$ & 1.99 & $2.73(-7)$ & 2.00  &   $1.77(-8)$ & 3.01   \\
$ 512\times 512\times 512$ &   6 & $1.71(-8)$ & 1.97   &  $6.90(-8)$ & 1.98 &   $2.18(-9)$ & 3.02    \\
\hline
  &  19.11 WU          & &  &  $8.68(-9)$ & &  $9.16(-10)$ &                           \\
\hline
\end{tabular}
\begin{tablenotes}
        \footnotesize
        \item[1] WU ({\it work unit}) is the computational cost of performing one relaxation sweep on the finest grid. Here, the $\textrm{EXCMG}_{bi-cg}$
        computation cost $= 6 + 48\times 2^{-3} + 384\times 2^{-6} + 533\times 2^{-9} + 285\times 2^{-12} \approx 19.11$.
      \end{tablenotes}
\end{threeparttable}
\label{table3}
\end{table}

~\\
{\bf Test Problem 4.}
The exact solution of the test problem 4 can be written as
\begin{equation}\label{test4}
u(x,y,z)= xyz \log(1+x+y+z).
\end{equation}

Applying the biharmonic operator on the exact solution, we can obtain the forcing term \emph{f}(\emph{x}, \emph{y}, \emph{z}) as follows:
\begin{equation}\label{force4}
\emph{f}(\emph{x}, \emph{y}, \emph{z})=\frac{-(2(4x^3 + 8x^2 + 15xyz + 4xy + 4xz + 4x + 4y^3 + 8y^2 + 4yz + 4y + 4z^3 + 8z^2 + 4z))}{(x + y + z + 1)^4}
\end{equation}

Obtain first boundary data from the exact solution while obtaining second boundary data \emph{$f_2$} by taking partial derivative for the exact solution.

Again, results listed on the table \ref{table4} of numerical experiments are performed with  EXCMG$_{bi-cg}$, using 3.6 thousand unknowns on the coarsest grid 32$\times$32$\times$32 and more than 135 million unknowns on the finest grid 512$\times$512$\times$512. Besides, from table \ref{table4}, we can see that numerical solution $u_h$ reaches almost full second-order accuracy, the initial guess $w_h$ is third-order approximation to numerical solution $u_h$, while the extrapolated solution $\tilde{u}_{h}$ increase the numerical solution's accuracy greatly. On the finest grid 512$\times$512$\times$512, using $w_h$ as numerical solution, the error ratio $r_h$ is already equal to 0.090. Thus we only need to perform iterations eight times to achieve the expected accuracy.

\begin{table}[!tbp]
\tabcolsep=6pt
\caption{Errors and convergence rates using EXCMG$_{bi-cg}$ for Problem 4.} \centering
\begin{threeparttable}
  \begin{tabular}{|c|c|cc|cc|cc|c|}
    \hline
    \multirow{2}{*}{Mesh}& \multirow{2}{*}{Iters} & \multicolumn{2}{c|}{$||u_h-u||_{2}$}  & \multicolumn{2}{c|}{$||u_h-u||_{\infty}$} & \multicolumn{2}{c|}{$||w_h-u_h||_{2}$}  \\
    \cline{3-8}&   &     Error &   Order  &     Error &   Order  &   Error &   Order  \\
\hline
$  32\times  32\times  32$ &   275 &$1.35(-6)$ &         &$3.47(-6)$ &       &  $2.19(-6)$ &        \\
$  64\times  64\times  64$ &   513 & $3.47(-7)$ & 1.96  &  $8.69(-7)$ & 2.00  &  $2.68(-7)$ & 3.03   \\
$ 128\times 128\times 128$ &   512 & $8.77(-8)$ & 1.98  &  $2.17(-7)$ & 2.00  &  $3.30(-8)$ & 3.02   \\
$ 256\times 256\times 256$ &   64 & $2.22(-8)$ & 1.98 & $5.46(-8)$ & 1.99  &   $4.11(-9)$ & 3.01   \\
$ 512\times 512\times 512$ &   8 & $5.70(-9)$ & 1.96   &  $1.39(-8)$ & 1.98 &   $5.11(-10)$ & 3.01    \\
\hline
  &  25.07 WU          & &  &  $1.50(-9)$ & &  $2.72(-10)$ &                           \\
\hline
\end{tabular}
\begin{tablenotes}
        \footnotesize
        \item[1] WU ({\it work unit}) is the computational cost of performing one relaxation sweep on the finest grid. Here, the $\textrm{EXCMG}_{bi-cg}$
        computation cost $= 8 + 64\times 2^{-3} + 512\times 2^{-6} + 513\times 2^{-9} + 275\times 2^{-12} \approx 25.07$.
      \end{tablenotes}
\end{threeparttable}
\label{table4}
\end{table}

~\\
{\bf Test Problem 5.}
The exact solution of the test problem 5 can be written as
\begin{equation}\label{test5}
u(x,y,z)= -e^{(10(x - 0.5)^2 + 10(y - 0.5)^2 + 10(z - 0.2)^2)}(- x^2 + x)(- y^2 + y)(- z^2 + z).
\end{equation}

Obtain first boundary data from the exact solution while obtaining second boundary data \emph{$f_2$} by taking partial derivative for the exact solution.

Again, we use five level grids which have 3.6 thousand unknowns on the coarsest grid 32$\times$32$\times$32 and more than 135 million unknowns on the finest grid 512$\times$512$\times$512. Additionally, from table \ref{table5}, we can see that numerical solution $u_h$ reaches almost full second-order accuracy, the initial guess $w_h$ is third-order approximation to numerical solution $u_h$, while the extrapolated solution $\tilde{u}_{h}$ increase the numerical solution's accuracy significantly.

\begin{table}[!tbp]
\tabcolsep=6pt
\caption{Errors and convergence rates using EXCMG$_{bi-cg}$ for Problem 5.} \centering
\begin{threeparttable}
  \begin{tabular}{|c|c|cc|cc|cc|c|}
    \hline
    \multirow{2}{*}{Mesh}& \multirow{2}{*}{Iters} & \multicolumn{2}{c|}{$||u_h-u||_{2}$}  & \multicolumn{2}{c|}{$||u_h-u||_{\infty}$} & \multicolumn{2}{c|}{$||w_h-u_h||_{2}$}  \\
    \cline{3-8}&   &     Error &   Order  &     Error &   Order  &   Error &   Order  \\
\hline
$  32\times  32\times  32$ &   432 &$8.86(-2)$ &         &$3.76(-1)$ &       &  $1.91(-1)$ &        \\
$  64\times  64\times  64$ &   873 & $2.42(-2)$ & 1.87  &  $1.01(-1)$ & 1.90  &  $2.55(-2)$ & 2.91   \\
$ 128\times 128\times 128$ &  1913 & $6.22(-3)$ & 1.96  &  $2.57(-2)$ & 1.97  &  $2.33(-3)$ & 3.45   \\
$ 256\times 256\times 256$ &   256 & $1.60(-3)$ & 1.96 & $6.45(-3)$ & 1.99  &   $2.16(-4)$ & 3.43   \\
$ 512\times 512\times 512$ &   32 & $4.45(-4)$ & 1.85   &  $1.65(-3)$ & 1.97 &   $2.45(-5)$ & 3.14    \\
\hline
  &  95.70 WU          & &  &  $2.84(-4)$ & &  $7.89(-5)$ &                           \\
\hline
\end{tabular}
\begin{tablenotes}
        \footnotesize
        \item[1] WU ({\it work unit}) is the computational cost of performing one relaxation sweep on the finest grid. Here, the $\textrm{EXCMG}_{bi-cg}$
        computation cost $= 32 + 256\times 2^{-3} + 1913\times 2^{-6} + 873\times 2^{-9} + 432\times 2^{-12} \approx 95.70$.
      \end{tablenotes}
\end{threeparttable}
\label{table5}
\end{table}

\section{Conclusion}

In this work, we propose a new extrapolation cascadic multigrid method $EXCMG_{bi-cg}$ to solve the linear three-dimensional biharmonic equation. By applying the Richardson extrapolation and quadratic interpolation methods on numerical solutions which are on current and previous grids, much better initial guesses of iterative solutions are obtained on the next finer grid so that the iterative time for Bi-CG solver is reduced. It is the main advantage of our work. Additionally, the introduction of the relative residual tolerance into our work enables us to obtain the desired accuracy conveniently. Furthermore, reducing computational time and the number of iteration, the numerical results of tests demonstrate that the $EXCMG_{bi-cg}$ method is efficient and particularly suitable for solving large scale problems.

\end{document}